\documentclass[12pt,a4paper,twoside,english]{amsart}

\usepackage{amsmath,amssymb,amsbsy,amsfonts,amsthm,latexsym,mathabx,
            amsopn,amstext,amsxtra,euscript,amscd,stmaryrd,mathrsfs,
            cite,array,mathtools,enumerate}
            
\usepackage[T1]{fontenc}
\usepackage[latin9]{inputenc}
\usepackage{babel}
\usepackage{verbatim}
\usepackage{amsthm}
\usepackage{amsmath}
\usepackage{amssymb}
 
 \usepackage[colorlinks,linkcolor=blue,anchorcolor=blue,citecolor=blue,backref=page]{hyperref}
\usepackage{color}

\usepackage{float}

\hypersetup{breaklinks=true}

\renewcommand*{\backref}[1]{}
\renewcommand*{\backrefalt}[4]{%
    \ifcase #1 (Not cited.)%
    \or        (p.\,#2)%
    \else      (pp.\,#2)%
    \fi}
    
    \def\({\left(}
\def\){\right)}

\def\mand{\qquad\mbox{and}\qquad}

\makeatletter

\pdfpageheight\paperheight 
\pdfpagewidth\paperwidth

\theoremstyle{plain}
\newtheorem{thm}{\protect\theoremname}[section]
 \newcommand\thmsname{\protect\theoremname}
 \newcommand\nm@thmtype{theorem}
 \theoremstyle{plain}

  \theoremstyle{remark}
  \newtheorem{rem}[thm]{\protect\remarkname}
  \theoremstyle{definition}
    \newtheorem{define}[thm]{\protect\definitionname}
  \newtheorem*{example*}{\protect\examplename}
  \theoremstyle{definition}
  \newtheorem{example}[thm]{\protect\examplename}
  \theoremstyle{plain}
  
  \theoremstyle{plain}
  
  \theoremstyle{plain}
  \newtheorem{cor}[thm]{\protect\corollaryname}

\makeatletter
  \theoremstyle{definition}
  \newtheorem{my@rem}[thm]{Remark}
  \renewenvironment{rem}{\begin{my@rem}}{\end{my@rem}}
\makeatother

\makeatother

  \providecommand{\examplename}{Example}
  \providecommand{\lemmaname}{Lemma}
  \providecommand{\propositionname}{Proposition}
  \providecommand{\remarkname}{Remark}
    \providecommand{\definitionname}{Definition}
  \providecommand{\theoremname}{Theorem}
\providecommand{\theoremname}{Theorem}
 \providecommand{\corollaryname}{Corollary}

\usepackage{mathtools}

\def\Q{{\mathbb Q}}
\def\R{{\mathbb R}}

\def\Z{{\mathbb Z}}
\def\G{{\mathbb G}}

\def\div{{\mathrm{div}}}
\def\Div{{\mathrm{Div}}}

\def\P{{\mathbb P}}
\def\C{{\mathbb C}}

\def\x{{\mathbf x}}
\def\a{{\mathbf a}}

\def\b{{\mathbf b}}

\newcommand{\h}{{\operatorname{h}}}

\def\Kc{K_{\mathrm c}}
\def\Kab{K_\mathrm{ab}}
\def\PhiX{\Phi_X}
\def\PhiXc{\Phi_{X,{\mathrm c}}}

\def\CVD{{\hfill\hfil{\lower 2pt\hbox{\vrule\vbox to 7pt
{\hrule width  5pt\varphifill\hrule}\varphirule}}}\par}

\def\Gm{\G_{\rm m}}

\begin{document}

\title[Abelian multiplicatively dependent points on curves]
{On abelian multiplicatively dependent points on a curve in a torus}

\author[A. Ostafe]{Alina Ostafe}
\address{A.O.: School of Mathematics and Statistics, University of New South Wales, Sydney NSW 2052, Australia}
\email{alina.ostafe@unsw.edu.au}

\author[M. Sha]{Min Sha}
\address{M.S.: Department of Computing, Macquarie University, Sydney, NSW 2109, Australia}
\email{shamin2010@gmail.com} 

\author[I.E. Shparlinski]{Igor E. Shparlinski}
\address{I.S.: School of Mathematics and Statistics, University of New South Wales, Sydney NSW 2052, Australia}
\email{igor.shparlinski@unsw.edu.au}

\author[U. Zannier]{Umberto Zannier}
\address{U.Z.: Scuola Normale Superiore, Piazza dei Cavalieri, 7, 56126 Pisa, Italy}
\email{u.zannier@sns.it}

\begin{abstract}  
We show, under some natural conditions,  that the set of  abelian (and thus also cyclotomic) multiplicatively dependent  points on
an irreducible curve over a number field is a finite union of preimages of roots of unity by a certain finite set of primitive characters from $\Gm^n$ to $\Gm$ 
restricted to the curve, and a finite set. 
We also introduce the notion of primitive multiplicative
dependence  and obtain a finiteness result for primitively multiplicatively dependent points
defined over a so-called Bogomolov extension of a number field. 
\end{abstract}

 


\maketitle

\section{Introduction} 

Let  $\Gm$ be  the multiplicative algebraic group over the complex numbers $\C$, that is $\Gm = \C^*$ endowed with the multiplicative group law. 
Let $n$ be a positive integer with $n \ge 2$. 
The points in $\Gm^n$ whose coordinates are roots of unity are called torsion points. 
In 1960s, Lang~\cite{Lang65}  conjectured that if a complex plane irreducible curve contains infinitely many torsion points, 
then the curve is a torsion coset of $\Gm^2$. 
This was soon confirmed by Ihara, Serre and Tate; see~\cite{Lang65} for more details. 

So far, there are two ways to generalise the above result. 
The first is to describe torsion points on an algebraic variety of higher dimension. 
This leads to the torsion points theorem proved independently by 
Laurent~\cite{Laurent} and Sarnak and Adams~\cite{SA}. 
The theorem asserts that the torsion points in a subvariety $Y \subseteq \Gm^n$ all lie and are Zariski dense in a finite number of torsion cosets contained in $Y$. 
Note that the torsion points constitute a multiplicative group of rank 0. 
So, the other way is to consider multiplicative groups of higher ranks. 
This was initiated by Bombieri, Masser and  Zannier~\cite{BMZ}. 
We refer to~\cite{Zan1} for more historic notes. 

More precisely, the paper~\cite{BMZ} studies the intersection of a geometrically irreducible  algebraic  curve $X \subseteq \Gm^n$, defined over a number field, and 
the union of proper algebraic subgroups of $\Gm^n$. 
As is well known (see, for example,~\cite[Corollary 3.2.15]{BG}),  each such subgroup $H$ is defined by a finite set of equations of the shape $x_1^{a_1}\cdots x_n^{a_n}=1$, with integer exponents not all zero. 
Hence, the above 
 intersection consists of the points $P=(\xi_1,\ldots ,\xi_n)\in X$ such that the coordinates $\xi_1=x_1(P),\ldots ,\xi_n=x_n(P)$ are all non-zero and are {\it multiplicatively dependent\/}.  In the sequel we  call such points  just {\it dependent}. 

 Let us note that any map  $(x_1,\ldots ,x_n)\mapsto x_1^{a_1}\cdots x_n^{a_n}$, denoted by $\x\to\x^\a$ in the sequel,  is a (rational)  homomorphism from $\Gm^n$ to $\Gm$, often called a {\it character}, which is non-trivial if and only if it is surjective, 
and in fact all rational  homomorphisms $\Gm^n\to\Gm$ are of this shape; see~\cite[Proposition 3.2.17]{BG}. The set of dependent points in $\Gm^n$ is just the union of the kernels of these non-trivial characters.

Coming back to $X$, one may assume for the problems in question  that it  is not contained in   any  proper algebraic subgroup of $\Gm^n$,  since otherwise  a (fixed) dependence occurs for all points of $X$. Under this assumption, $X\cap H$ is finite for any such algebraic subgroup $H$, so that  it  consists of  points  defined over $\overline\Q$. Also, on varying $H$ it is easy to see that the union of these intersections, that is, the set of dependent points on $X$,  is infinite. 

It has been proved in~\cite[Theorem~1]{BMZ},   that, under the assumption that {\it $X$ is not contained in
any translate of a proper algebraic subgroup of $\Gm^n$}, the absolute logarithmic {\it Weil\/}  height in $X\cap H$ is bounded independently of $H$; see~\cite{BG,Zan2} for a
background on heights.

This   more stringent  assumption means that  no  equation 
$x_1^{a_1}\cdots x_n^{a_n}=c$, for a constant $c$ and integers $a_i$ not all zero,   can hold identically on $X$, that is, $X$ is not contained in any fibre of a non-trivial rational character to $\Gm$.
 In~\cite{BMZ}, this hypothesis has been noted to be  necessary for the bounded-height conclusion. 
 Throughout the paper, we
always  assume this condition for $X$. 

Now, one can ask what can be further said on restricting the dependent  points of $X$ to be defined over proper subfields of $\overline \Q$. For instance, over any number field, or even just imposing  some bound for the degree over $\Q$,  the mentioned  bounded-height result of~\cite{BMZ} implies that there can be only finitely many such points on   $X$ which are dependent,  because of the {\it Northcott theorem\/}, see~\cite[Theorem~1.6.8]{BG}. 

Other fields relevant to this context are the {\it maximal cyclotomic fields\/} $\Kc$,  that is, 
fields obtained by adjoining all roots of unity to a number field $K$.  This is
especially for the reason that  any point   with a coordinate which is a root of unity is automatically dependent.  
The goal of this paper is, roughly,   to describe the structure of abelian (and thus also cyclotomic) points on $X$ which are  dependent. 
Note for instance  that   the said result of~\cite{BMZ} does not lead us to expect  automatically finiteness, since the Northcott Property certainly fails for the maximal cyclotomic fields.

 In fact,  a situation which generates infinitely many  dependent  cyclotomic  points on $X$ is as follows. 
Suppose that the character  $\varphi:\x=(x_1, \ldots ,x_n)\to \x^\a$,
where $\x^\a = x_1^{a_1}\cdots x_n^{a_n}$, restricts on $X$  to a birational correspondence $\varphi_X:X\to \Gm$.  
Then, the inverse image $\varphi_X^{-1}(\zeta)$ of any root of unity $\zeta$  ``usually'' consists of  one point  $P\in X$, hence necessarily defined over a cyclotomic extension, and certainly by definition this point is  dependent, since  $\varphi(P)^m=1$ for some $m>0$.   More generally, if  $\varphi_X$ is a factor of an isogeny $\varphi_X\circ\rho: \Gm\to\Gm$ (with $\rho$ some birational map $\Gm\to X$ defined over $\Kc$), again the inverse image of any root of unity  produces   cyclotomic points, automatically dependent. 

More precisely, here we obtain the following result:
 In Theorem~\ref{T.main}, we 
prove that the above situation  characterises, in finite terms, all of the cyclotomic dependent points in $X$, 
and in fact we obtain a more general result for {\it abelian\/} dependent points in $X$, 
that is, for points from the abelian closure $\Kab$ of $K$. 
 In Theorem~\ref{thm:XKc}, we  obtain a similar result for  a more stringent relation 
however for an even broader class of fields, called fields with the  Bogomolov propery.

For possible applications of our results see Remark~\ref{rem:Appl} below and~\cite{OSSZ}.

\section{Main results} 

To state in precision our  results, let  $K$ be a number field and let $X\subseteq \Gm^n$ be a geometrically irreducible curve defined over $K$. Let us denote by $U\subset\overline\Q$ the set of all roots of unity. 
As before we use  $\Kc$ to denote the maximal cyclotomic extension of $K$, that is, $\Kc=K(U)$ and
use $\Kab$ to denote the maximal abelian extension over $K$ 
which includes but is generally much larger than $\Kc$ (if $K\neq \Q$).
 
 We  also denote by $\varphi_X: X\to \Gm$ 
 the restriction to $X$ of the character $\varphi$ on $\Gm^n$.  Recall that a character is called {\it primitive} if it is not of the shape $\psi^m$ for a character $\psi$ and an integer $m>1$. 

As usual, for any field $L$ with $K \subseteq L$ we write $X(L)$ 
for the set of $L$-rational points on $X$. 

Throughout the paper, the  ``height'' means  the ``absolute logarithmic Weil height'', 
which  we denote by $\h:\overline\Q\to \R\cap[0, \infty)$; 
see~\cite{BG,Zan2}. 

Let $\PhiX$ be the set of   primitive characters $\varphi:\Gm^n\to\Gm$ with the 
property that there exists a birational isomorphism $\rho:\Gm\to X$ such that   $\varphi_X\circ \rho$ is an isogeny of $\Gm$.  

We also use $\PhiXc$  to denote the subset of $\PhiX$ consisting of those characters
for which $\rho$ can be defined over $\Kc$

We now present our main result:

\begin{thm} \label{T.main}  
Suppose that $X$ is not contained in any translate of a proper algebraic subgroup of $\Gm^n$. 
Then, $\PhiX$ is a  finite set, and  the set of dependent points in $X(\Kab)$ is the union of  
the set $\bigcup_{\varphi\in\PhiXc}\varphi_X^{-1}(U)$  
and a finite set.
\end{thm}

\begin{rem}[\textit{Cyclotomic dependent points}]
\label{rem:CDP}
Since  such birational isomorphisms $\rho$  corresponding to characters in the set $\PhiXc$ are defined over $\Kc$, it is easy to see that 
for any $\varphi \in \PhiXc$, we have $\varphi_X^{-1}(U) \subseteq X(\Kc)$. 
So, the set of dependent points in $X(\Kc)$ is also the union of  
 the set $\bigcup_{\varphi\in\PhiXc}\varphi_X^{-1}(U)$  and a finite set. 
\end{rem}

Taking into account Remark~\ref{rem:CDP} we immediately derive:

\begin{cor} \label{cor:Xab Xc}  Under the conditions of Theorem~\ref{T.main}, 
 the set of  dependent points in $X(\Kab)$ is the union of a finite set with the set  of dependent  points in $X(\Kc)$.
\end{cor}

\begin{rem}[{{\it Curves of positive genus}}]
\label{rem:g>0}
 Note that if   $\PhiXc$ is  empty,  Theorem~\ref{T.main}  implies that the set of dependent   points in $X(\Kab)$  is   finite. 
 For instance, this automatically occurs if  $X$ has genus $g_X>0$. 
 Indeed, to see this, it is enough to notice that $\Gm$ cannot be birationally isomorphic to $X$ when $g_X>0$.  
 However, even if  $g_X = 0$, for $\PhiXc$ to be nonempty we need the severe condition that some monomial in the coordinate functions on $X$ is a power of a rational function on $X$ of degree $1$. 
\end{rem}

\begin{rem}[{{\it Elliptic curves with complex multiplication}}]
\label{rem:elliptic}
Let $E$ be an elliptic curve defined by an affine Weierstrass equation over $K$ which has complex multiplication. 
Without loss of generality, we can assume that $K$ contains the endomorphism ring of $E$. 
Then, it is well-known that (see, for instance,~\cite{Serre}) the field generated over $K$ by all the $E$-torsion points (that is, torsion points in the sense of elliptic curves) is contained in $\Kab$. 
So, from Remark~\ref{rem:g>0}, we know that there are only finitely many $E$-torsion points which are also dependent. 
\end{rem}

\begin{rem}[{{\it Tightness of the assumptions}}]
\label{rem:Cond}
The assumption of Theorem~\ref{T.main}  that $X$ is not contained in any translate of a proper algebraic subgroup cannot be  replaced with the weaker (obviously necessary) one that $X$ is not contained in any proper algebraic subgroup.  This is shown for instance by the example  
$$
X=\{(x_1,x_2,x_3)\in\Gm^3~:~x_1=2,\  x_3^2=x_2^3+1\}.
$$ 
Then, $X$ has genus $1$ (so $\PhiXc$ is empty), whereas  the infinitely many   points 
$(2,2^m,\sqrt{2^{3m}+1})\in X(\Q_{\mathrm c})$, $m\in\Z$,    
are clearly dependent; also, it is very easy to check that $X$ is not contained in any proper algebraic subgroup. 
Here, $\sqrt{2^{3m}+1} \in \Q_{\mathrm c}$ due to the Kronecker-Weber theorem. 
\end{rem}

\begin{rem} [{{\it Effectiveness}}]
The proof  in Section~\ref{sec:proofs} below shows that if $X$ is given effectively,  then $\PhiXc$ and the finite set appearing in Theorem~\ref{T.main} may be effectively computed. 
\end{rem}

\begin{rem} [{{\it Nontriviality of isogenies}}]
In general, the isogeny mentioned in the definition of $\PhiXc$ cannot be taken  trivial. 
This is equivalent to saying that $\PhiXc$ cannot be generally replaced with the set of characters whose restriction to $X$ induces  a birational isomorphism $X\cong\Gm$. For an explicit case, 
see the following  Example~\ref{exam:x x-1}.
\end{rem}

\begin{example}
\label{exam:x x-1}
 Let $X$ be defined in $\Gm^2$ by $x_1=(x_2-1)^d$. One can check that $\PhiXc$ consists of the  characters $x_1^{\pm 1},x_2^{\pm 1}$. Note that setting $x_1=\zeta$, a root of unity, indeed yields the cyclotomic dependent points $(\zeta,\mu+1)$ on $X$, where $\mu^d=\zeta$. Similarly, we may set $x_2=\zeta$ and get $((\zeta-1)^d,\zeta)$. For any number field $K$, Theorem~\ref{T.main} implies that these families 
of points account for all but finitely many dependent points in $X(\Kc)$. 
In particular, we recover that equations of the shape $x^m(x-1)^n=1$
 cannot ``generally'' be solved within $\Kc$. 
\end{example}

One  can also supplement  Theorem~\ref{T.main} with even more explicit descriptions of the dependent points, which would follow as corollaries of its statement (and partially of its proof). 
In turn, this may lead to other consequences. Here, we do not pursue in this task and  limit to the natural  question  whether   a dependence relation $\xi_1^{a_1}\cdots \xi_n^{a_n}=1$ for $P=(\xi_1,\ldots ,\xi_n)\in X(\overline{\Q})$ can be {\it primitive}, that is, with coprime 
exponents $\gcd(a_1,\ldots,a_n)=1$. 
We say that the corresponding point $P$ is {\it primitively dependent\/}.   

In fact, the set of primitively dependent points on $X$ is nothing else than the intersection of $X$ with the \textit{union of connected proper algebraic subgroups of $\Gm^n$}.

We need to recall the notion of a field with the  {Bogomolov Property\/}; see~\cite{ADZ}:

\begin{define} 
We say that a subfield $L$ of $\overline\Q$ has the {\it Bogomolov Property} if
there exists a   constant 
$C(L)>0$ which  depends only on  $L$, 
such that for any $\alpha \in L^*\setminus U$ we have $\h(\alpha) \ge C(L)$.
\end{define}

For example, in view of~\cite{AZ}, we can choose $L=\Kc$, or even $L=\Kab$, see~\cite{ADZ} for further examples. 

Then we have the following: 

\begin{thm} \label{thm:XKc} 
Under  the assumptions of Theorem~\ref{T.main} on $X$,  
if $L$ is a field with the Bogomolov Property over $K$, then
there are only finitely many primitively dependent points in $X(L)$.
\end{thm}

The same example as in   Remark~\ref{rem:Cond} shows that again the assumptions on $X$ cannot be weakened. 
Also, it is easy to see that the set of primitively dependent points in $X(\overline\Q)$ is always infinite. 

As an analogue of Remark~\ref{rem:elliptic}, we record the following remark. 

\begin{rem} [{{\it Elliptic curves}}]
Let $E$ be an elliptic curve defined by an affine Weierstrass equation over $\Q$. 
By a result of Habegger~\cite[Theorem~1]{Hab2}, the field generated over $\Q$ by all the $E$-torsion points satisfies the Bogomolov Property. 
So, it follows from Theorem~\ref{thm:XKc} that there are only finitely many $E$-torsion points which are also primitively dependent. 
\end{rem}

\section{Proofs}
\label{sec:proofs}

\subsection{Notation} 
We also recall that the notation $A \ll B$  (sometimes we write this also as $B \gg A$) 
is equivalent to the inequality 
$|A| \leq c B$ for some constant $c$, which throughout the paper may depend on $K$ and $X$.

\subsection{Proof of Theorem~\ref{T.main}} 

We divide the proof into three steps. 
We first show that $\PhiX$ is finite, which is our {\it Step~(I)}, and then we prove the second assertion.
From Remark~\ref{rem:CDP}, we know that  the set $\bigcup_{\varphi\in\PhiXc}\varphi_X^{-1}(U)$  is contained in $X(\Kab)$ 
and its elements are dependent points. 
Besides, any dependent point in $X(\Kab)$ is, by definition, mapped to a root of unity under some character $\varphi$ of $\Gm^n$.  
So, for the second assertion, we only need to prove that there are only finitely many choices of such a character $\varphi$, which we do at {\it Step~(II)}, and 
if there are infinitely many points  in $X(\Kab)$ sent to a root of unity by $\varphi$,   then $\varphi \in \PhiXc$, which is done at {\it Step~(III)}.

{\it Step~(I)}: Proving  the finiteness of $\PhiX$. 

Let $\widetilde X$ be a smooth projective curve defined over $K$ and $K$-birational to $X$.  For a character $\varphi$, the restriction $\varphi_X$ to $X$ is a rational function on $X$, so it is also  a rational function on $\widetilde X$ and we may consider its divisor
$\div(\varphi_X)$ in the set  of divisors $\Div(\widetilde X)$ of $\widetilde X$; 
see~\cite[Appendix~A]{BG} for a background on divisors. 

We then have a map (a homomorphism) $\varphi\mapsto \div(\varphi_X)\in \Div(\widetilde X)$ 
which associates to a (rational)  character  $\varphi:\Gm^n\to\Gm$ the divisor in $\widetilde X$ of   $\varphi_X$.  
Note that this map is injective, since any $\varphi$ in the kernel yields a constant $\varphi_X$, which cannot happen  in virtue of our assumption on $X$. 
Then, it suffices to show that the image of this map on $\PhiX$ is finite.

Let  $\varphi\in\PhiX$, so we may write $\varphi_X\circ\rho(t)=t^d$, where $\rho:\Gm\to X$ is a birational isomorphism and $t$ is a coordinate on $\G_{m}$.  
In fact, we may extend $\rho$ to  an isomorphism $\rho:\P_1\to\widetilde X$, and the same equation holds on viewing $\varphi_X$ as a map 
from $\widetilde X$ to the projective line $\P_1$. 
Hence, it appears that $\varphi_X$ can have only one zero and one pole, and so its divisor is   of the shape $m((P)-(Q))$ for $P,Q$ distinct points on $\widetilde X$. But  since $\varphi$ is a monomial in the coordinates, these $P,Q$ necessarily lie among the zeros and poles of the coordinate functions $x_1,\ldots ,x_n$ viewed as functions on $\widetilde X$, so $P,Q$ have only finitely many possibilities. 

Finally, for given $P,Q$, let $\varphi,\psi\in\PhiX$ correspond to $P,Q$, such that 
$$
\div\(\varphi_X\)=\ell ((P)-(Q)) \mand \div\(\psi_X\)=m((P)-(Q))
$$ 
for integers $\ell$ and $m$.  Then $\div((\varphi^m/\psi^\ell)_X)=0$, so in fact $\varphi^m=\psi^\ell$ by the above injectivity. 
But both $\varphi,\psi$ are primitive and this yields $\ell=\pm m$. 
Hence, each pair $P,Q$ can give rise to at most two elements of $\PhiX$, and, 
combined with already established finiteness of the choices for $P$ and $Q$, 
the finiteness of  $\PhiX$ follows.

{\it Step~(II)}: Describing the characters.

To go ahead, we start with some arguments similar to those 
for~\cite[Theorem~2]{BMZ}  (see also~\cite[Chapter~2]{Zan1}).

Let $P\in X(\Kab)$ be a dependent point, so the coordinates $\xi_i=x_i(P)$ generate a multiplicative subgroup  $\Gamma_P\subset  (\Kab)^*$ of rank $r\le n-1$. We may thus write (on invoking elementary abelian group decomposition) 
\begin{equation}
\label{eq:coord}
\xi_i=\zeta_i\prod_{j=1}^rg_j^{m_{ij}}, \qquad  i =1, \ldots, n,
\end{equation} 
for  generators $g_i\in\Q(P)$ of the
 torsion-free part of $\Gamma_P$,  integers $m_{ij}$ and roots of unity $\zeta_i\in\Q(P)$.   
 
Now, by an argument coming from the geometry of numbers we may actually  find the generators $g_i$ so that, for any integers $b_1,\ldots ,b_r$, 
\begin{equation}\label{E.lemma2}
\h(g_1^{b_1}\cdots g_r^{b_r})\ge c_r \(|b_1|\h(g_1)+\ldots +|b_r|\h(g_r)\),
\end{equation}
where $c_r>0$ is a positive number depending only on $r$; see~\cite[Lemma~2]{BMZ} (which in turn refers to a previous  result of Schlickewei~\cite{Schl}). 

Now, by the already mentioned~\cite[Theorem~1]{BMZ}, the height of dependent points in $X(\overline\Q)$ is uniformly  bounded, so that 
$$
\h\(\prod_{j=1}^rg_j^{m_{ij}}\)\ll 1.
$$ 
Using then  the inequality~\eqref{E.lemma2} to bound below the height of $\xi_i$ in~\eqref{eq:coord}, we find
\begin{equation}\label{E.ub}
\sum_{j=1}^r|m_{ij}|\h(g_j)\ll 1,
\end{equation}
where the implicit constant depends only on $X$ and the ambient dimension $n$. 

 In~\cite{BMZ}, this kind of inequality is exploited with $r\le n-2$ on using certain lower bounds for heights coming from a higher-dimensional version of lower bounds of Dobrowoski, due to Amoroso and David~\cite{AD}.  Here we may well have $r=n-1$, however we may take advantage of the fact that the $g_i$ are in $\Kab$ (and certainly they are not roots of unity since they are multiplicatively independent). For these fields an absolute lower bound for the height has been 
 proved by Amoroso and Zannier~\cite{AZ}. Specifically, by~\cite[Theorem~1.2]{AZ}, applied  with  
  $\alpha=g_j$  we have
 \begin{equation}\label{E.AZ} 
  \h(g_j)\gg 1,
 \end{equation}
 where now the implicit constant depends only on $K$ (in fact only on the degree $[K:\Q]$). 
 So, here we indeed need that the point $P$ is defined over $\Kab$ 
 (the other place we need this is in the end of {\it Step~(III)} below). 
 
 Combining~\eqref{E.ub} with~\eqref{E.AZ}, we obtain that the absolute values $|m_{ij}|$ are upper bounded independently of any dependent point $P\in X(\Kab)$.  
 We may now consider separately the finitely many possibilities which arise for the $m_{ij}$ and thus  in the sequel we may suppose that the $m_{ij}$ are fixed, that is, independent of $P$.  
 
 Let now $(b_1,\ldots ,b_n)\in\Z^n$ be a nonzero integer vector orthogonal to the $(m_{1j},\ldots ,m_{nj})$, $j=1,\ldots, r$. Since $r<n$, 
  such an integer vector exists, and we may view it as 
fixed, like for the integers $m_{ij}$; we may also take it primitive (that is, with coprime coordinates). 
 The previous equations~\eqref{eq:coord} for $\xi_1, \ldots \xi_n$ yield
 \begin{equation}\label{E.ker}
 \xi_1^{b_1}\cdots \xi_n^{b_n}=\zeta_P,
 \end{equation}
 where $\zeta_P$ is a root of unity (which  depends on $P$).   
  Let us denote by 
  $$
  \pi :\x\mapsto \x^\b
  $$ 
  the character so obtained, which is primitive in our above meaning. 
  
 Note that this character $\pi$ is taken from a prescribed finite set, and it  sends the point $P$ to a root of unity; in fact,  this  is already an approach to Theorem~\ref{T.main}. 
  What is missing are the properties of the character stated in  Theorem~\ref{T.main}.  
 So, to complete the proof, we only need to show that if infinitely many points 
 $P\in X(\Kab)$ are sent to a root of unity by $\pi$,    
then $\pi$ is in $\PhiXc$.

{\it Step~(III)}: Proving that $\pi \in \PhiXc$.

 For this we use~\cite[Theorem~2.1]{Zan0} (we note that in the case of $X(\Kc)$ the 
 result of~\cite[Theorem~1]{DZ} is also sufficient). 
 
 Let $T_k$ be the set of torsion points of $\Gm^k$, $k\ge 1$. 
 For a rational map $\tau$ from a geometrically irreducible variety $Y$ to $\Gm^k$, 
 the (PB) condition in~\cite{Zan0} means that for any integer $m > 0$, the pullback $\Gm^k \times_{[m],\tau} Y$ is geometrically irreducible, 
 where $[m]$ is the $m$-th power map. 
 The result in~\cite[Theorem~2.1]{Zan0} asserts that:
 \begin{equation}  \label{eq:PB}
 \begin{array}{l}
 \textit{if $\tau$ is a cover $($that is, dominant rational map of finite degree$)$} \\
 \textit{defined over $\Kc$ and satisfies the {\rm (PB)} condition, there exists a}  \\
 \textit{finite union $W$ of proper torsion cosets such that if $v \in T_k \setminus W$,}\\
  \textit{then $v \in \tau(Y)$ and if $\tau(u)=v$, then $[\Kc(u) : \Kc] = \deg \tau$.} 
 \end{array}
 \end{equation}

 So, in order to apply~\eqref{eq:PB}, we need to construct a suitable cover $\tau$ related to the character $\pi$. 
 For this, we first factor $\pi_X$ according to~\cite[Proposition~2.1]{Zan0} into an isogeny of algebraic groups and a rational map satisfying the (PB) condition, 
 from which we can construct a suitable cover. Then, we choose a point $u$ such that 
 we can compute the degree  $[\Kc(u) : \Kc]$ from two different ways 
 (including~\eqref{eq:PB}) and yield that $\deg \tau =1$ and so $\tau$ is an isomorphism, 
 which leads to $\pi \in \PhiXc$. 
 
 We first show that $\pi_X$ is a cover. 
In the equation~\eqref{E.ker}, fixing a root of unity for the right-hand side gives a proper algebraic subgroup $H$ of $\Gm^n$. 
Since $X\cap H$ is finite and there are infinitely many $P\in X(\Kab)$ such that $\pi_X(P)$ is a root of unity, 
we must have that $U \cap \pi_X(X(\Kab))$ is infinite, and thus it is Zariski dense in $\Gm$ 
because obviously (by the finiteness of the roots of a nontrivial 
univariate polynomial) any closed subset in  $\Gm$ is either 
finite or is $\Gm$ itself.
Hence, $\pi_X: X \to \Gm$ is a cover.  

We can now factor $\pi_X$ as $\lambda_0 \circ \rho_0$ according to the 
second claim 
of~\cite[Proposition~2.1]{Zan0} such that 
$\lambda_0: Y \to \Gm$ is an isogeny of algebraic groups and $\rho_0: X \to Y$ is a rational map satisfying the (PB) condition. 
Note that here $(X,Y, \Gm)$ correspond to $(Y,Z,X)$ in the notation of~\cite[Proposition~2.1]{Zan0}. 
 Since there is a dual isogeny $\widehat{\lambda_0}: \Gm \to Y$ of $\lambda_0$, we see that 
 $Y$ is isomorphic to a quotient of $\Gm$ by a finite subgroup. 
 Now note that a finite subgroup of 
 $\Gm$  is in fact a subgroup of roots of unity. So, it is a kernel of an $m$-th power 
 map $[m]: \Gm \to \Gm$, 
that is,
\begin{equation}    \label{eq:PowMap}
 [m]: x \mapsto x^m,
\end{equation}
for some integer $m \ge 1$.
Under the power map, we can see that $\Gm$  modulo the kernel (that is, the finite subgroup) 
is isomorphic to $\Gm$. Therefore, $Y$ is in fact  isomorphic to $\Gm$ as algebraic groups.
Hence, we in fact can factor 
\begin{equation}    \label{eq:factor}
\pi_X = [m] \circ \rho,
\end{equation}
where $\rho: X \to \Gm$ is a rational map satisfying the (PB) condition, and  the 
 $m$-th power map  $[m]: \Gm \to \Gm$ is given by~\eqref{eq:PowMap}
(see~\cite[Proposition~3.2.17]{BG}). 
Since there are only finitely many such characters $\pi$ (and also $\pi_X$), and so there are only finitely many such integers $m$, 
we add all $m$-th roots of unity into $K$. In other words, we enlarge the field $K$ by adding finitely many roots of unity. 
Since $\pi$ and $[m]$ are defined over $\Q$, by~\eqref{eq:factor} we know that $\rho$ is defined over $K$ (alternatively, over $\Kc$ if we do not enlarge the field $K$).

Furthermore, for any point $P \in X(\Kab)$ and any Galois automorphism $\sigma$ over $K$, denoting by $P^\sigma$ the image of $P$ under $\sigma$, 
by definition we have 
\begin{equation}   \label{eq:sigma}
\rho(P^\sigma) = \rho(P)^\sigma.
\end{equation}

Since $\pi_X$ is a cover,  the map $\rho$ is also a cover. 
Indeed, if  $\rho$ is not a cover, then  by definition the image of $\rho$ is not dense in  $\Gm$, which
by~\eqref{eq:factor} implies that 
the image of $\pi_X$ is   not dense either, contradicting to the fact that  $\pi_X$ is a cover. 
Now, we define the product cover 
$$
\psi=\rho \times \rho: X \times X \to \Gm^2 \cong \Gm \times \Gm, 
$$
which is also defined over $K$.   
Since the cover $\rho$ satisfies the (PB) condition, the same is true for the cover $\psi$. 
Clearly, we have $\deg \psi = (\deg \rho)^2$. 
The cover $\psi$ is exactly what we want for invoking~\eqref{eq:PB}. 
The rest is to choose a suitable point. 

Since $X$ is defined over $K$, if $P \in X(\Kab)$, then $P^\sigma \in X(\Kab)$ for any Galois automorphism $\sigma$ over $K$. 
Note that for any point $P \in X(\Kab)$, by~\eqref{eq:factor}, $\pi_X(P)$ is a root of unity if and only if  $\rho(P)$ is a root of unity. 
So, in view of the infinity of $U \cap \pi_X(X(\Kab))$, we know that $U \cap \rho(X(\Kab))$ is also infinite. 
Thus, considering the set $S$ of images of points of the form $(P,P^\sigma)$ under $\psi$ for any $P \in \rho^{-1}(U) \cap X(\Kab)$ and any Galois automorphism $\sigma$ over $K$, 
$S$ is an infinite set and is a torsion subset of $\Gm^2$, whose elements by~\eqref{eq:sigma} have the form  $(\zeta,\zeta^\sigma)$ for some root of unity $\zeta$. 
Applying~\eqref{eq:PB} to both $\rho$ and $\psi$ and noticing that each proper torsion coset of $\Gm^2$ contains only finitely many torsion points of such form  $(\zeta,\zeta^\sigma)$,  
there exist an element $(\zeta,\zeta^\sigma) \in S$ and a point $(P,P^\sigma)$ with $P \in X(\Kab)$ mapped to $(\zeta,\zeta^\sigma)$ under $\psi$ such that 
$$
[K(P):K] =\deg \rho,  \qquad [K(P,P^\sigma) : K] = \deg \psi = (\deg \rho)^2. 
$$
On the other hand, since $P \in X(\Kab)$, the field extension $K(P)/K$ is automatically normal, whence it  contains $K(P^\sigma)$ for any conjugate $P^\sigma$ of $P$ 
over $K$. So, we obtain 
$$
[K(P,P^\sigma) : K] = [K(P) : K] = \deg \rho. 
$$
Hence, we must have $\deg \rho =1$. 
That is, $\rho$ is an isomorphism. 
So, by~\eqref{eq:factor}, we have $\pi_X \circ \rho^{-1} = [m]$, which is an isogeny of $\Gm$, 
and thus $\pi \in \PhiXc$. 
This completes the proof.

 \subsection{Proof of Theorem~\ref{thm:XKc}}
 
 Let $P=(\xi_1,\ldots ,\xi_n)\in X(L)$ be a primitively dependent point, and let $\xi_1^{a_1}\cdots \xi_n^{a_n}=1$ be a primitive relation, that is, with the $a_i$ coprime integers.

We could indeed use the result of Theorem~\ref{T.main} for the proof 
 if $L=\Kab$, but it is as easy to go again through the proof. Arguing exactly as in the previous proof 
  (the equation~\eqref{E.AZ} still holds since $L$ has the Bogomolov Property), we find bounded integers $m_{ij}$, $1\le i\le n$, $1\le j\le r$,  and generators $g_j$ for the group $\xi_1^\Z\cdots \xi_n^\Z$  as above. Here $r$ is again the rank of such group. It follows that the vectors $(m_{1j},\ldots ,m_{nj})$, $j=1,\ldots, r$, are linearly independent.
  
  Note that $(a_1,\ldots ,a_n)$ must be orthogonal to the $(m_{1j},\ldots ,m_{nj})$, $j=1,\ldots, r$, by the independence of the $g_j$. 
 
 If $r=n-1$, then  the space orthogonal to the $(m_{1j},\ldots ,m_{nj})$, $j=1,\ldots, r$, has dimension $1$, hence $(a_1,\ldots ,a_n)$ is uniquely determined up to sign, since it is primitive. Hence, since the $m_{ij}$ are bounded, it has only finitely many possibilities, and the equation $x_1^{a_1}\cdots  x_n^{a_n}=1$ yields only finitely many points $P\in X$.
 
 Therefore we may assume that $r\le n-2$. Now, we could conclude immediately by using~\cite[Theorem~2]{BMZ}, which yields finiteness for points verifying two or more independent multiplicative relations. But we can conclude in an easier way, without using such result. 
 
 For this, let us start  by observing that a standard argument yields two independent vectors 
 $$
 (b_1,\ldots ,b_n) \mand (c_1,\ldots ,c_n)
 $$
  orthogonal to the vectors $(m_{1j},\ldots ,m_{nj})$, $j=1,\ldots, r$, and bounded independently of $P$ (since the $m_{ij}$ are likewise bounded). 
  Hence, we may assume that these vectors are fixed for an infinity of the points $P$ in question. 
 
 Let $\varphi,\psi$ be the corresponding characters, which are multiplicatively independent. Then, $\varphi(P)$ and $\psi(P)$ are both roots of unity. The map $(\varphi,\psi):\Gm^n\to\Gm^2$ is a homomorphism sending the curve $X$ to a curve $Y\subset \Gm^2$. 
 The curve $Y$ then contains the torsion points $(\varphi(P),\psi(P))$. If these points make up an infinite set, then $Y$ is a translate of a torus. But then, by the multiplicative independence of $\varphi$ and $\psi$,  $X$ is contained in a nontrivial translate of an algebraic subgroup, against the assumptions. 
 Hence, there are only finitely many points $\varphi(P)$, which then implies that $P$ lies in a finite set.

\section{Possible generalisations and applications}

 \begin{rem}[{{\it Higher dimensions}}] 
 One may ask what happens for higher dimensional varieties in place of the curve $X$. Now the bounded-height result of~\cite{BMZ}, 
 which is a crucial tool for the present proofs,  is not true in the most obvious generalization, but a correct analogue has been proved by Habegger~\cite{Hab}. 
 Probably this leads to higher-dimensional analogues of the present conclusions. 
 \end{rem}

\begin{rem}[{{\it Multiplicative dependence of rational values}}]
\label{rem:Appl}
We recall that any $n$ rational functions 
over $K$  in one variable $t$, not all constant, always parametrize a curve $X$, which is rational over $K$. Hence, our results can be 
used to investigate multiplicative dependence of the values on $\Kc$ (or more generally on $\Kab$) of given rational functions in $K(t)$. To satisfy the necessary condition of Theorem~\ref{T.main}, we 
have to make a natural assumption that there is no nontrivial product  
of these rational functions which is a constant function. 
We refer to~\cite{OSSZ} for some partial results. 
\end{rem}

\section*{Acknowledgement}

The authors would like to thank the referee for giving constructive comments which helped to improve the paper. 
During the preparation of this work, A. O. was supported by the UNSW FRG Grant PS43704, M. S. was supported by the Macquarie University Research Fellowship,  
and I. S.   was  supported   by the ARC Grant DP170100786.

\end{document}